\documentclass[11pt]{article}

\usepackage{bbm}
\begin{document}

\title{Modular invariance,  characteristic numbers and $\eta$ invariants}
\author{Fei Han\  and\  Weiping Zhang}
\date{}
\maketitle

\begin{abstract}
We generalize the ``miraculous cancellation" formulas of
Alvarez-Gaum\'e, Witten and  Kefeng Liu to a twisted version
where an extra complex line bundle is involved. We also apply our
result  to discuss intrinsic relations between the higher
dimensional Rokhlin type congruence formulas of Ochanine [O1],
Finashin [F] and Zhang [Z1, 2]. In particular, an analytic proof
of the Finashin congruence formula is given.
\end{abstract}

$$\ $$

\noindent{\bf \S 1. Introduction}

$\ $

A well-known theorem of Rokhlin [R1] states that the Signature of
an oriented closed smooth spin 4-manifold is divisible by 16. It
has two kinds of higher dimensional generalizations. One is due
to Atiyah and Hirzebruch [AtH] stating that the
$\widehat{A}$-genus of an $8k+4$ dimensional oriented closed
smooth spin manifold is an even integer. The other one, due to
Ochanine (cf. [O1]), states that the Signature of such an $8k+4$
dimensional  manifold is divisible by 16.

Recall that Atiyah and Hirzebruch's result in [AtH] goes beyond
the $\widehat{A}$-genus. In fact, they  proved that if $E$ is a
real vector bundle over an $8k+4$ dimensional oriented closed
smooth spin manifold $M$, then $\langle \widehat{A}(TM){\rm
ch}(E\otimes {\mathbbm C}),[M]\rangle\in 2{\mathbbm Z}$.

It turns out that these two generalizations of the original
Rokhlin divisibility are closely related: in [L], Landweber shows
 how one can use the ideas of elliptic genus to deduce  the
Ochanine divisibility directly from the divisibility results of
Atiyah and Hirzebruch.

Now let $M$ be a closed oriented smooth 4-manifold not necessarily
spin. Let $B$ be an orientable characteristic submanifold of $M$,
that is, $B$ is a compact two dimensional submanifold of $M$ such
that $[B]\in H_2(M,{\mathbbm Z}_2)$ is dual to the second
Stiefel-Whitney class of $TM$.

In this case, Rokhlin established in [R2] a congruence formula of
the type
$${{\rm Sign}(M)-{\rm Sign}(B\cdot B)\over 8}\equiv \phi (B)\ \ \ {\rm
mod}\ \ 2{\mathbbm Z},\eqno (1.1)$$ where $B\cdot B$ is the
self-intersection of $B$ in $M$ and $\phi(B)$ is a spin cobordism
invariant associated to $(M,B)$.

Clearly, when $M$ is spin and $B=\emptyset$, (1.1) reduces to the
original Rokhlin divisibility.

In the case where $B$ might be nonorientable, an extension of
(1.1) was proved by Guillou and Marin in [GuM] (see  [KT] for a
comprehensive  account).

In [O1], Ochanine generalizes (1.1) to $8k+4$ dimensional closed
Spin$^c$ manifolds. His formula is also of type (1.1) and thus
extends his divisibility result to Spin$^c$ manifolds.

In [F], Finashin generalizes the Ochanine congruence to the case
where the characteristic submanifold $B$ is nonorientable. His
formula also extends the Guillou-Marin congruence to any $8k+4$
dimensional closed oriented manifold.

 On the other hand, Zhang [Z1, 2] proves another type of
congruence formulas for $8k+4$ dimensional  manifolds. His
results generalize the Atiyah-Hirzebruch divisibilities and, when
restricted to 4-manifolds, also recover (1.1) as well as the
Guillou-Marin extension of it.

In view of Landweber's elliptic genus proof of the Ochanine
divisibility, it is natural to ask whether there exist similar
intimate relations between these higher dimensional
generalizations of  (1.1) as well as its extension by Guillou and
Marin.

In [Z2, Appendix],  this question was briefly dealt with   for
the case of 12 dimensional Spin$^c$ manifolds, with the help of
the so called ``miraculous cancellation" formula first proved by
the physicists  Alvarez-Gaum\'e and Witten [AGW] (cf. (2.43) in
the text). In particular, an intrinsic analytic interpretation of
Ochanine's spin cobordism invariant $\phi (B)$
is given. Moreover, it was pointed  out that a higher dimensional
generalization of the ``miraculous cancellation" formula of
Alvarez-Gaum\'e and Witten would lead to a better understanding
of the general congruence formulas due to Ochanine [O1] and
Finashin [F].

This required   higher dimensional ``miraculous cancellation"
formula  was established  by Kefeng Liu [Li] by developing modular
invariance properties of characteristic forms. In some sense,
Liu's formula refines the  argument of Landweber [L] to the level
of differential forms.

In [LiZ], by combining  Liu's formula with the analytic arguments
in [Z2], Liu and Zhang were able to give an intrinsic analytic
interpretation of the Ochanine invariant $\phi(B)$ for any $8k+2$
dimensional closed spin manifold $B$ (Compare with [O2]), as well
as the Finashin invariant [F] for any $8k+2$ dimensional closed
pin$^-$ manifold appearing in the Finashin congruence formula. In
particular, this leads to  analytic versions, stated in [LiZ,
Theorems 4.1, 4.2], of the Finashin and Ochanine congruences.

Now the remained  question is whether these  analytic versions of
the Finashin and Ochanine congruences could be deduced directly
from the Rokhlin type congruences in [Z1, 2].

The purpose of this paper is to give a positive answer to this
question. To be more precise, what we get is a generalization of
the ``miraculous cancellation" formulas of Alvarez-Gaum\'e, Witten
and Liu to a twisted version where an extra complex line bundle
is involved. Modular invariance properties developed in [Li]
still play an important role in the proof of such an extended
cancellation formula. When applying our formula to Spin$^c$
manifolds, we are led directly to an unexpected refined version
of [LiZ, Theorem 4.2] (cf. (3.2) in the text). Moreover, by
combining this twisted cancellation formula with the analytic
Rokhlin type congruences proved in [Z2], we are able to give a
direct analytic proof of [LiZ, Theorem 4.1], which by [LiZ,
Theorem 3.1] is equivalent to the original Finashin congruence
formula.

The rest of the article is organized as follows. In Section 2 we
establish our twisted extension of the ``miraculous cancellation"
formulas of Alvarez-Gaum\'e, Witten and Liu. We include also an
Appendix to this section where we state an analogous  twisted
cancellation formula for $8k$ dimensional manifolds.
In Section 3, we apply the twisted cancellation formula proved in
Section 2 to $8k+4$ dimensional Spin$^c$ manifolds and show how
it leads directly to a refinement of the Ochanine congruence
formula. Finally, in Section 4, we combine the results in Section
2 with the Rokhlin type congruences of Zhang [Z2, Theorem 3.2] to
give a direct proof of the analytic version of the Finashin
congruence stated in [LiZ, Theorem 4.1].

Parts of the results of this article have been announced in [HZ].

$\ $

\noindent {\bf Acknowledgment} This work was partially supported
by MOEC and the 973 project.

$$\ $$

\noindent {\bf \S 2. Modular invariance and a twisted ``miraculous
cancellation" formula}

$\ $

In this section, we generalize the ``miraculous cancellation"
formulas of Alvarez-Gaum\'e, Witten and  Kefeng Liu ([AGW], [Li])
to a twisted version where an extra complex line bundle (or,
equivalently, a rank two real oriented  vector bundle) is
involved.

This section is organized as follows. In a), we present the basic
geometric data and recall the definitions of the characteristic
forms to be discussed. In b), we state the main result of this
section, which is a twisted extension of the ``miraculous
cancellation" formulas of Alvarez-Gaum\'e, Witten and Liu. In c),
we recall some basic facts about the modular forms which will be
used in d) to give a proof of the main result stated in b).
Finally, in e), we specialize the main result stated in b) to the
tangent bundle case and give an explicit expression of it in the
12 dimensional case. There is also an appendix to this section
where we include an analogous twisted cancellation formula for
$8k$ dimensional manifolds.

$\ $

{\bf a). Some characteristic forms}

$\ $

Let $M$ be an $8k+4$ dimensional Riemannian manifold. Let
$\nabla^{TM}$ be the associated Levi-Civita connection and
$R^{TM}=\nabla^{TM,2}$ the curvature of $\nabla^{TM}$.

Let $\widehat{A}(TM, \nabla^{TM})$, $\widehat{L}(TM, \nabla^{TM})$
be the Hirzebruch characteristic forms defined by
$${ \widehat{A}(TM, \nabla^{TM})
={\det}^{1/2}\left({{\sqrt{-1}\over  4\pi}R^{TM} \over
\sinh\left({ \sqrt{-1}\over 4\pi}R^{TM}\right)}\right),\ \
\widehat{L}(TM, \nabla^{TM}) ={\det}^{1/2}\left({{\sqrt{-1}\over
2\pi}R^{TM} \over \tanh\left({ \sqrt{-1}\over
4\pi}R^{TM}\right)}\right)}. \eqno (2.1)$$

Let $E$, $F$ be two Hermitian vector bundles over $M$ carrying
with Hermitian connections $\nabla^E$, $\nabla^F$ respectively.
Let $R^E=\nabla^{E,2}$ (resp. $R^F=\nabla^{F,2}$) be the
curvature of $\nabla^E$ (resp. $\nabla^F$). If we set the formal
difference $G=E-F$, then  $G$ carries with an induced Hermitian
connection $\nabla^G$ in an obvious sense. We define the
associated Chern character form  as
$${\rm ch}(G,\nabla^G)={\rm tr}\left[\exp\left({\sqrt{-1}\over 2\pi}R^E\right)\right]
-{\rm tr}\left[\exp\left({\sqrt{-1}\over
2\pi}R^F\right)\right].\eqno(2.2)$$

In the rest of this paper, when there will be no confusion about
the Hermitian connection $\nabla^E$ on a Hermitian vector bundle
$E$, we will also write simply ${\rm ch}(E)$ for the associated
Chern character form.

For any complex number $t$, let $$\Lambda_t(E)={\mathbbm
C}|_M+tE+t^2\Lambda^2(E)+\cdots ,\ \ \ S_t(E)={\mathbbm
C}|_M+tE+t^2S^2(E)+\cdots$$  denote respectively the total
exterior and symmetric powers  of $E$.

We recall  the following relations between these two operations
(cf. [At, Chap. 3]),
$$ S_t(E)=\frac{1}{\Lambda_{-t}(E)},\ \ \ \
 \Lambda_t(E-F)=\frac{\Lambda_t(E)}{\Lambda_t(F)}.\eqno(2.3)$$
Therefore,  we have the following formulas for Chern character
forms,
 $$ {\rm ch}(S_t(E) )=\frac{1}{{\rm ch}(\Lambda_{-t}(E) )},\ \ \ \
{\rm ch}(\Lambda_t(E-F) )=\frac{{\rm ch}(\Lambda_t(E) )}{{\rm
ch}(\Lambda_t(F) )}\ .\eqno(2.4)$$

If $W$ is a  real Euclidean vector bundle over $M$ carrying with a
Euclidean connection $\nabla^W$,   then its complexification
$W_\mathbbm{C}=W\otimes \mathbbm{C}$ is a complex vector bundle
over $M$ carrying with a canonically induced Hermitian metric from
the Euclidean metric of $W$, as well as a Hermitian connection
from $\nabla^W$.

We refer to [Z3, Sect. 1.6] for the definitions and notations of
other Pontrjagin (resp. Chern) forms associated to real (resp.
complex) vector bundles with connections.

$\ $

{\bf b). A twisted ``miraculous cancellation" formula}

$\ $

We make the same assumptions and use the same notations as in a).

Let $V$ be a rank $2l$ real Euclidean vector bundle over $M$
carrying with a Euclidean connection $\nabla^V$.

Let $\xi$ be a rank two real oriented Euclidean vector bundle
over $M$ carrying with a Euclidean connection $\nabla^\xi$.

If $E$ is a complex vector bundle over $M$, set
$\widetilde{E}=E-{\mathbbm C}^{{\rm rk}({E})}$.

Let $q=e^{2\pi \sqrt{-1}\tau}$ with $\tau\in {\mathbbm H}$, the
upper half complex plane.

Set $$\Theta_1(T_{\mathbbm C}M,V_{\mathbbm C},\xi_{\mathbbm
C})=\bigotimes_{n=1}^\infty S_{q^n}(\widetilde{T_{\mathbbm C}M})
\otimes \bigotimes_{m=1}^\infty
\Lambda_{q^m}(\widetilde{V}_{\mathbbm
C}-2\widetilde{\xi}_{\mathbbm C})\otimes
\bigotimes_{r=1}^\infty\Lambda_{q^{r-{1\over
2}}}(\widetilde{\xi}_{\mathbbm
C})\otimes\bigotimes_{s=1}^\infty\Lambda_{-q^{s-{1\over
2}}}(\widetilde{\xi}_{\mathbbm C}),$$
$$\Theta_2(T_{\mathbbm C}M,V_{\mathbbm C},\xi_{\mathbbm C})=\bigotimes_{n=1}^\infty S_{q^n}(\widetilde{T_{\mathbbm C}M})
\otimes \bigotimes_{m=1}^\infty \Lambda_{-q^{m-{1\over
2}}}(\widetilde{V}_{\mathbbm C}-2\widetilde{\xi}_{\mathbbm C})
\otimes \bigotimes_{r=1}^\infty\Lambda_{q^{r-{1\over
2}}}(\widetilde{\xi}_{\mathbbm
C})\otimes\bigotimes_{s=1}^\infty\Lambda_{q^{s}}(\widetilde{\xi}_{\mathbbm
C}).\eqno(2.5)$$ Clearly,  $\Theta_1(T_{\mathbbm C}M,V_{\mathbbm
C},\xi_{\mathbbm C }) $ and $\Theta_2(T_{\mathbbm C}M,V_{\mathbbm
C},\xi_{\mathbbm C})$ admit formal Fourier expansion  in
$q^{1/2}$ as
$$\Theta_1(T_{\mathbbm C}M,V_{\mathbbm C},\xi_{\mathbbm C})=A_0(T_{\mathbbm C}M,V_{\mathbbm C},\xi_{\mathbbm C})+
A_1(T_{\mathbbm C}M,V_{\mathbbm C},\xi_{\mathbbm
C})q^{1/2}+\cdots,$$
$$\Theta_2(T_{\mathbbm C}M,V_{\mathbbm C},\xi_{\mathbbm C})=B_0(T_{\mathbbm C}M,V_{\mathbbm C},\xi_{\mathbbm C})
+B_1(T_{\mathbbm C}M,V_{\mathbbm C},\xi_{\mathbbm
C})q^{1/2}+\cdots,\eqno(2.6)$$ where the $A_j$'s and $B_j$'s are
elements in the semi-group formally generated by Hermitian vector
bundles over $M$. Moreover, they carry canonically induced
Hermitian connections.

Let $c=e(\xi,\nabla^\xi)$ be the  Euler form of $\xi$ canonically
associated to $(\xi,\nabla^\xi)$ (cf. [Z3, Sect. 3.4]). Let
$R^V=\nabla^{V,2}$ denote the curvature of $\nabla^V$.

If $\omega$ is a differential form over $M$, we denote by
$\omega^{(8k+4)}$ its top  degree component.

We can now state our main result of this section as follows.

$\ $

\noindent{\bf Theorem 2.1} {\it If the equality for the first
Pontrjagin forms $p_1(TM,\nabla^{TM})=p_1(V,\nabla^V)$ holds, then
one has the  equation for $(8k+4) $-forms,
$$\left\{ {\widehat{A}(TM,\nabla^{TM}){\det}^{1/2}
\left(2\cosh\left({\sqrt{-1}\over 4\pi}R^V\right)\right) \over
\cosh^2({c\over 2})}\right\}^{(8k+4)} $$ $$=2^{l+2k+1}\sum_{r=0}^k
2^{-6r}\left\{ \widehat{A}(TM,\nabla^{TM}){\rm ch}(
b_r(T_{\mathbbm C}M,V_{\mathbbm C},\xi_{\mathbbm C })
)\cosh\left({c\over 2}\right)\right\}^{(8k+4)},\eqno(2.7)
$$ where each $b_r(T_{\mathbbm C}M,V_{\mathbbm C},\xi_{\mathbbm C })
$, $0\leq r\leq k$, is a  canonical integral linear combination
of  $B_j(T_{\mathbbm C}M,V_{\mathbbm C},\xi_{\mathbbm C }) $,
$0\leq j\leq r$.}

$\ $

Certainly, when $\xi={\bf R}^2$ and $c=0$, Theorem 2.1 is exactly
Liu's result in [Li, Theorem 1], which, in the
$(TM,\nabla^{TM})=(V,\nabla^V)$ and $k=1$ case, recovers the
original ``miraculous cancellation" formula of Alvarez-Gaum\'e
and Witten [AGW].

Theorem 2.1 will be proved in d) by using  Liu's argument in [Li],
in taking account  into the appearance of $\xi$ and $c$. In the
next subsection, for the sake of self completion, we will recall
some of the materials concerning modular forms which will be used
in d).

$\ $

{\bf c). Some properties about the Jacobi theta functions and
modular forms}

$\ $

We first recall that the four Jacobi theta functions are defined
as follows (cf. [C]):
$$\theta(v,\tau)=2q^{1/8}\sin(\pi v)
\prod_{j=1}^\infty\left[(1-q^j)(1-e^{2\pi
\sqrt{-1}v}q^j)(1-e^{-2\pi \sqrt{-1}v}q^j)\right]\ ,\eqno(2.8)$$
 $$\theta_1(v,\tau)=2q^{1/8}\cos(\pi v)
 \prod_{j=1}^\infty\left[(1-q^j)(1+e^{2\pi \sqrt{-1}v}q^j)
 (1+e^{-2\pi \sqrt{-1}v}q^j)\right]\ ,\eqno(2.9)$$
 $$\theta_2(v,\tau)=\prod_{j=1}^\infty\left[(1-q^j)
 (1-e^{2\pi \sqrt{-1}v}q^{j-1/2})(1-e^{-2\pi \sqrt{-1}v}q^{j-1/2})\right]\ ,\eqno(2.10)$$
$$ \theta_3(v,\tau)=\prod_{j=1}^\infty\left[(1-q^j)
(1+e^{2\pi \sqrt{-1}v}q^{j-1/2})(1+e^{-2\pi
\sqrt{-1}v}q^{j-1/2})\right]\ ,\eqno(2.11)$$ where $q=e^{2\pi
\sqrt{-1}\tau}$ with $\tau\in {\mathbbm H}$.

Set
$${\theta}'(0,\tau)=\left. {\partial \theta (v,\tau)\over \partial
v}\right|_{v=0}.\eqno(2.12)$$

We refer to [C, Chap. 3] for a proof of the following Jacobi
identity.

$\ $

\noindent {\bf Proposition 2.2} {\it The following identity holds,
$${\theta}'(0,\tau)= \pi\,
\theta_1(0,\tau)\theta_2(0,\tau)\theta_3(0,\tau)\, .\eqno(2.13)$$}

Let as usual $SL_2(\mathbbm{Z})= \{(\begin{array}{cc}a&b\\
c&d\end{array})\, |\ a,b,c,d\in\mathbbm{Z},\ ad-bc=1\}$. Let $S=(\begin{array}{cc}0&-1\\
1&0\end{array})$, $T=(\begin{array}{cc}1&1\\
0&1\end{array})$ be the two generators of $SL_2(\mathbbm{Z})$.
They act on ${\mathbbm H}$ by $S\tau =-1/\tau$, $T\tau=\tau+1$.

One has the following transformation laws of theta functions under
$S$ and $T$ (cf. [C]),
$$\theta(v,\tau+1)=\theta(v,\tau),\ \ \
\theta\left(v,-{1}/{\tau}\right)=\sqrt{-1}\tau^{1/2}e^{-\tau
v^2}\theta\left(\tau v,\tau\right)\ ;\eqno(2.14)$$
$$\theta_1(v,\tau+1)=\theta_1(v,\tau),\ \ \
\theta_1\left(v,-{1}/{\tau}\right)=\tau^{1/2}e^{-\tau
v^2}\theta_2(\tau v,\tau)\ ;\eqno(2.15)$$
$$\theta_2(v,\tau+1)=\theta_3(v,\tau),\ \ \
\theta_2\left(v,-{1}/{\tau}\right)=\tau^{1/2}e^{-\tau
v^2}\theta_1(\tau v,\tau)\ ;\eqno(2.16)$$
$$\theta_3(v,\tau+1)=\theta_2(v,\tau),\ \ \
\theta_3\left(v,-{1}/{\tau}\right)=\tau^{1/2}e^{-\tau
v^2}\theta_3(\tau v,\tau)\ .\eqno(2.17)$$

Let $\Gamma $ be a subgroup of $SL_2(\mathbbm{Z})$.

$\ $

\noindent {\bf Definition 2.3} A modular form over $\Gamma$ is a
holomorphic function $f(\tau)$ on ${\mathbbm H}$ such that
$$f(g\tau)
:=f\left(\frac{a\tau+b}{c\tau+d}\right)
=\chi(g)(c\tau+d)^kf(\tau)\ \ \ \  \mbox{for\  any}\ \
g\in \left(\begin{array}{cc}a&b\\
c&d\end{array}\right)\in \Gamma,\eqno(2.18)$$
 where $\chi:\Gamma\rightarrow\mathbbm{C}^*$ is a character of
 $\Gamma$ and $k$ is called the weight of $f$.

 $\ $

Let ${\mathcal M}_{\mathbbm R}(\Gamma)$ denote the ring of
modular forms over $\Gamma$ with real Fourier coefficients.

Following [Li], denote by $\theta_j=\theta_j(0,\tau)$, $1\leq
j\leq 3$, and define
$$\delta_1(\tau)={1\over
8}(\theta_2^4+\theta_3^4),\ \ \ \ \varepsilon_1(\tau)={1\over
16}\theta_2^4\theta_3^4,\eqno(2.19)$$ $$\delta_2(\tau)=-{1\over
8}(\theta_1^4+\theta_3^4),\ \ \ \ \varepsilon_2(\tau)={1\over
16}\theta_1^4\theta_3^4.\eqno(2.20)$$
 They admit Fourier expansion
$$\delta_1(\tau)={1\over 4}+6q+\cdots,\ \ \ \ \varepsilon_1(\tau)={1\over
16}-q+\cdots,\eqno(2.21)$$
$$\delta_2(\tau)=-{1\over 8}-3q^{1/2}+\cdots,\ \ \ \
\varepsilon_2(\tau)=q^{1/2}+\cdots,\eqno(2.22)$$ where the
``$\cdots$" terms are higher degree terms all having integral
coefficients. They also verify the following transformation laws
under $S$ (cf. [L] and [Li]),
$$\delta_2(-{1/ \tau})
=\tau^2\delta_1(\tau),\ \ \ \ \varepsilon_2(-{1/
\tau})=\tau^4\varepsilon_1(\tau).\eqno(2.23)$$

Let $\Gamma_0(2)$, $\Gamma^0(2)$ be two subgroups of
$SL_2(\mathbbm{Z})$ defined by
$$ \Gamma_0(2)=\left\{\left.\left(\begin{array}{cc}
a&b\\
c&d
\end{array}\right)\in SL_2(\mathbbm{Z})\,\right|\,c\equiv 0\ \ {\rm mod} \ \ 2{\mathbbm Z}\right\},$$

$$ \Gamma^0(2)=\left\{\left.\left(\begin{array}{cc}
a&b\\
c&d
\end{array}\right)\in SL_2(\mathbbm{Z})\,\right|\,b\equiv 0\ \ {\rm mod} \ \ 2{\mathbbm Z}\right\}.$$
Then $T,\ ST^2ST$ are the two generators of $\Gamma_0(2)$, while
$STS,\ T^2STS$ are the two generators of $\Gamma^0(2)$.

The following weaker version of [Li, Lemma 2] will be used in the
next subsection.

$\ $

\noindent {\bf Lemma 2.4} {\it One has that $\delta_2$ (resp.
$\varepsilon_2$) is a modular form of weight $2$ (resp. $4$) over
$\Gamma^0(2)$. Furthermore, one has ${\mathcal M}_{\mathbbm
R}(\Gamma^0(2) )={\mathbbm
R}[\delta_2(\tau),\varepsilon_2(\tau)]$.}

$\ $

{\bf d). A proof of Theorem 2.1}

$\ $

Without loss of generality, we will adopt the Chern roots
formalism as in [Li] in the computation of characteristic forms.

Recall that if $\{ w_i \}$ are the formal Chern roots of a
Hermitian vector bundle $E$ carrying with a Hermitian connection
$\nabla^E$, then one has the following formula for the Chern
character form of the exterior power of $E$
(Compare with [Hi]),
$${\rm ch}(\Lambda_t(E) )=\prod_i(1+e^{w_i}t).\eqno(2.24)$$

Set for $\tau\in {\mathbbm H}$ and $q=e^{2\pi \sqrt{-1}\tau}$ that
$$P_1(\tau)=\left\{ {\widehat{A}(TM,\nabla^{TM}){\det}^{1/2}
\left(2\cosh\left({\sqrt{-1}\over 4\pi}R^V\right)\right) \over
\cosh^2({c\over 2})} {\rm ch}\left(\Theta_1(T_{\mathbbm
C}M,V_{\mathbbm C },\xi_{\mathbbm C
}),\nabla^{\Theta_1(T_{\mathbbm C}M,V_{\mathbbm C },\xi_{\mathbbm
C })}\right)\right\}^{(8k+4)},\eqno(2.25)$$
$$P_2(\tau)=\left\{ \widehat{A}(TM,\nabla^{TM}){\rm
ch}\left(\Theta_2(T_{\mathbbm C}M,V_{\mathbbm C },\xi_{\mathbbm C
}),\nabla^{\Theta_2(T_{\mathbbm C}M,V_{\mathbbm C},\xi_{\mathbbm C
})}\right)\cosh\left(c\over
2\right)\right\}^{(8k+4)},\eqno(2.26)$$ where
$\nabla^{\Theta_i(T_{\mathbbm C}M,V_{\mathbbm C},\xi_{\mathbbm
C})}$, $i=1,\ 2$, are the induced
 Hermitian connections with $q^{j/2}$-coefficients on $\Theta_i(T_{\mathbbm C}M,V_{\mathbbm C},\xi_{\mathbbm C})$
 from those on the ${A_j}$'s and ${B_j}$'s (Compare with  (2.6)).

Since (2.7) is a local formula over $M$, without loss of
generality, we may well assume that both $TM$ and $V$ are
oriented. Let $\{\pm 2\pi \sqrt{-1}y_v\}$ (resp. $\{\pm2\pi
\sqrt{-1}x_j\}$) be the formal Chern roots for $(V_{\mathbbm{C}},
\nabla^{V_{\mathbbm{C}}})$ (resp. $(
T_{\mathbbm{C}}M,\nabla^{T_{\mathbbm{C}}M})$). Let $c=2\pi
\sqrt{-1}u$.

{}From (2.1) and (2.25), one finds,
$$P_1(\tau)=2^l\left\{ \left(\prod_{j=1}^{4k+2}\frac{\pi
x_j}{\sin (\pi x_j)}\right) \left(\prod_{v=1}^l
 \cos (\pi y_v)\right)\frac{{\rm ch}
 \left(\Theta_1(T_{\mathbbm C}M,V_{\mathbbm C},\xi_{\mathbbm C})\right)}{\cos^2\left(\pi u\right)}
 \right
 \}^{(8k+4)}\, .
 \eqno(2.27)$$

{}From (2.4) and (2.5), one can write  ${\rm ch}
 (\Theta_1(T_{\mathbbm C}M,V_{\mathbbm C},\xi_{\mathbbm C}))$ as follows,
\newcommand{\h}{\begin{eqnarray*}}
 \newcommand{\e}{\end{eqnarray*}}
 \newcommand{\n}{{\rm ch}(\Lambda_{q^n}}
 \newcommand{\m}{{\rm ch}(\Lambda_{q^m}}
 \newcommand{\rr}{{\rm ch}(\Lambda_{q^{r-\frac{1}{2}}}}
 \newcommand{\s}{{\rm ch}(\Lambda_{-q^{s-\frac{1}{2}}}}
 \newcommand{\nn}{\prod_{n=1}^\infty}
 \newcommand{\mm}{\prod_{m=1}^\infty}
 \newcommand{\rrr}{\prod_{r=1}^\infty}
 \newcommand{\sss}{\prod_{s=1}^\infty}
 \newcommand{\jj}{\prod_{j=1}^{4k+2}}
 \newcommand{\vv}{\prod_{v=1}^l}
$${\rm ch}(\Theta_1(T_{\mathbbm C}M,V_{\mathbbm C},\xi_{\mathbbm
C}) ) = \nn \frac{{\rm ch}(\Lambda_{-q^n} (\mathbbm{C}^{8k+4})
)}{{\rm ch}(\Lambda_{-q^n} (T_\mathbbm{C}M) )} \mm \frac{{\rm
ch}(\Lambda_{q^m} (V_\mathbbm{C}) )}{\m(\mathbbm{C}^{2l}) )}
\prod_{t=1}^\infty\left(\frac{{\rm
ch}(\Lambda_{q^t}(\mathbbm{C}^2) )}{{\rm
ch}(\Lambda_{q^t}(\xi_\mathbbm{C}) )}\right)^2 $$ $$\rrr
\frac{\rr(\xi_\mathbbm{C}) )}{\rr (\mathbbm{C}^2) )}\sss
\frac{\s(\xi_\mathbbm{C}) )}{\s (\mathbbm{C}^2))}\, .
 \eqno(2.28)$$
{}From (2.8), the Jacobi identity (2.13) and (2.24), one deduces
directly that
\newcommand{\0}{\theta}
\newcommand{\1}{\theta_1}
\newcommand{\2}{\theta_2}
\newcommand{\3}{\theta_3}
$$\prod_{j=1}^{4k+2}\frac{\pi
x_j}{\sin (\pi x_j)}\nn \frac{{\rm ch}(\Lambda_{-q^n}
(\mathbbm{C}^{8k+4}) )}{{\rm ch}(\Lambda_{-q^n} (T_\mathbbm{C}M)
)}=\prod_{j=1}^{4k+2}x_j \frac{\pi
\1(0,\tau)\2(0,\tau)\3(0,\tau)}{\0(x_j,\tau)} =\jj x_j
\frac{\0'(0,\tau)}{\0(x_j,\tau)}\, .\eqno(2.29)$$ Similarly, from
(2.9)-(2.11) and (2.24), one deduces that
$$\vv \cos (\pi y_v) \mm \frac{\m (V_\mathbbm{{C}}) )}{\m (
\mathbbm{C}^{2l} ) )} =\vv \frac{\1(y_v,\tau)}{\1(0,\tau)},\ \ \ \
\rrr \frac{\rr(\xi_\mathbbm{C}) )}{\rr (\mathbbm{C}^2)
)}=\frac{\3(u,\tau)}{\3(0,\tau)}\, ,\eqno (2.30)$$
$$\frac{1}{\cos^2\left(\pi u\right)}
 \prod_{t=1}^\infty\left(\frac{{\rm ch}(\Lambda_{q^t}(\mathbbm{C}^2)
 )}{{\rm ch}(\Lambda_{q^t}(\xi_\mathbbm{C}) )}\right)^2=
 \frac{\1^2(0,\tau)}{\1^2(u,\tau)},\ \ \ \
 \sss
\frac{\s(\xi_\mathbbm{C}))}{\s
(\mathbbm{C}^2))}=\frac{\2(u,\tau)}{\2(0,\tau)}\, .\eqno (2.31)$$

Putting (2.27)-(2.31) together, one finds that the first part of
the following result holds.

$\ $

\noindent {\bf Proposition 2.5} {\it The following two identities
hold,
$$P_1(\tau)=2^l \left\{ \prod_{j=1}^{4k+2}\left(x_j{\theta'(0,\tau)\over \theta(x_j,\tau)}\right)
 \left(\prod_{v=1}^{l}{\theta_1(y_v,\tau)\over \theta_1(0,\tau)}\right) {\theta_1^2(0,\tau)\over \theta_1^2(u,\tau)}
 {\theta_3(u,\tau)\over \theta_3(0,\tau)}
{\theta_2(u,\tau)\over
\theta_2(0,\tau)}\right\}^{(8k+4)},\eqno(2.32)$$
$$P_2(\tau)= \left\{ \prod_{j=1}^{4k+2}\left(x_j{\theta'(0,\tau)\over \theta(x_j,\tau)}\right)
 \left(\prod_{v=1}^{l}{\theta_2(y_v,\tau)\over \theta_2(0,\tau)} \right){\theta_2^2(0,\tau)\over \theta_2^2(u,\tau)}
 {\theta_3(u,\tau)\over \theta_3(0,\tau)}
{\theta_1(u,\tau)\over
\theta_1(0,\tau)}\right\}^{(8k+4)}.\eqno(2.33)$$ }

{\it Proof}. Formula (2.32) has been proved above. By doing
similar computation, one also gets (2.33).  Q.E.D.

$\ $

 Next, by direct verifications in applying  the transformation laws
(2.14)-(2.17) to (2.32), (2.33) respectively, one gets

$\ $

\noindent {\bf Proposition 2.6} {\it If the equation for the first
Pontrjagin forms $p_1(TM,\nabla^{TM})=p_1(V,\nabla^V)$ holds, then
$P_1(\tau)$ is a modular form of weight $4k+2$ over
$\Gamma_0(2)$; while $P_2(\tau)$ is a modular form of weight
$4k+2$ over $\Gamma^0(2)$. Moreover, the following identity holds,
$$P_1\left({-{1/ \tau}}\right) =2^l\tau^{4k+2}P_2(\tau).\eqno(2.34)$$}

We can now proceed to prove Theorem 2.1 as follows.

Observe that at any point $x\in M$, up to the volume form
determined by the metric on $T_xM$, both $P_i(\tau)$, $i=1,\ 2$,
can be viewed as a power series of $q^{1/2}$ with real Fourier
coefficients. Thus, one can combine Lemma 2.4 and Proposition 2.6
to get, at $x$,  that
$$P_2(\tau)=h_0(8\delta_2)^{2k+1}+h_1
(8\delta_2)^{2k-1}\varepsilon_2+\cdots+h_k(8\delta_2)\varepsilon_2^k\
, \eqno(2.35)$$ where each $h_j$, $0\leq j\leq k$, is a real
multiple of the volume form at $x$.

By (2.23),  (2.34) and (2.35), one deduces that
$$P_1(\tau)=\frac{2^l}{\tau^{4k+2}}P_2\left(-{1/ \tau}\right)
=\frac{2^l}{\tau^{4k+2}}
\left[h_0\left(8\delta_2\left(-{1/\tau}\right)\right)^{2k+1}
+h_1\left(8\delta_2\left(-{1/\tau}\right)\right)^{2k-1}\varepsilon_2\left(-{1/\tau}\right)
\right.$$
$$\left. +\cdots
+h_k\left(8\delta_2\left(-{1/\tau}\right)\right)\left(\varepsilon_2\left(-{1/\tau}\right)\right)^k\right]
=2^l\left[h_0(8\delta_1)^{2k+1}+h_1(8\delta_1)^{2k-1}\varepsilon_1+\cdots+h_k(8\delta_1)\varepsilon_1^k\right]\,
 .\eqno(2.36)$$

{}By (2.5), (2.21), (2.25) and by setting $q=0$ in (2.36), one
deduces that $$\left\{ {\widehat{A}(TM,\nabla^{TM}){\det}^{1/2}
\left(2\cosh\left({\sqrt{-1}\over 4\pi}R^V\right)\right) \over
\cosh^2\left({c\over 2}\right)}\right\}^{(8k+4)}=
2^{l+2k+1}\sum_{r=0}^k2^{-6r}h_r\, .\eqno(2.37)$$

Now in order to prove (2.7), one need to show that each $h_r$,
$0\leq r\leq k$, can be expressed through a canonical integral
linear combination of  $\{ \widehat{A}(TM,\nabla^{TM}){\rm
ch}(B_j,\nabla^{B_j}) \cosh ({c\over 2})\}^{(8k+4)}$, $0\leq
j\leq r$, with coefficients not depending on $x\in M$.

As in [L], one can use the induction method to prove this fact
easily by comparing the coefficients of $q^{j/2}$, $j\geq 0$,
between the two sides of (2.35). We leave the details to the
interested reader.

Here for convenience  we write out the explicit expressions for
$h_0$ and $h_1$ as follows.
$$h_0=-\left\{\widehat{A}(TM,\nabla^{TM}) \cosh \left({c\over 2}\right)
\right\}^{(8k+4)},\eqno(2.38)$$
$$h_1=\left\{\widehat{A}(TM,\nabla^{TM})\left[24(2k+1)-{\rm
ch}(B_1,\nabla^{B_1})\right] \cosh \left({c\over
2}\right)\right\}^{(8k+4)}.\eqno(2.39)$$

$\ $

 \noindent {\bf Remark 2.7}  From (2.6), (2.22), (2.26),
(2.35) and Theorem 2.1, one finds that for any integer $r\geq 0$,
$\{ \widehat{A}(TM,\nabla^{TM}){\rm ch}(B_r,\nabla^{B_r}) \cosh
({c\over 2})\}^{(8k+4)}$ can be expressed through a canonical
integral linear combination of  $\{
\widehat{A}(TM,\nabla^{TM}){\rm ch}(B_j,\nabla^{B_j}) \cosh
({c\over 2})\}^{(8k+4)}$, $0\leq j\leq k$. This fact is by no
means trivial for $r\geq k+1$. It depends heavily on the modular
invariance of $P_2(\tau)$.

$\ $

{\bf e). The case of $V=TM$}

$\ $

In this subsection we apply Theorem 2.1 to the case where $V=TM$
and $\nabla^V=\nabla^{TM}$. In this case, in view of (2.1), (2.7)
becomes
$$\left\{ {\widehat{L}(TM, \nabla^{TM})\over \cosh^2({c\over 2})}\right\}^{(8k+4)}
=8\sum_{r=0}^k2^{6k-6r}\left\{ \widehat{A}(TM,\nabla^{TM}){\rm
ch}\left( b_r(T_{\mathbbm C}M,\xi_{\mathbbm C }) \right)
\cosh\left({c\over 2}\right)\right\}^{(8k+4)},\eqno(2.40)$$ where
we have used the simplified notation $b_r(T_{\mathbbm
C}M,\xi_{\mathbbm C })$ for $ b_r(T_{\mathbbm C}M,T_{\mathbbm
C}M,\xi_{\mathbbm C })$.

Now we assume $k=1$, that is, $\dim M=12$. Then by concentrating
on the coefficients of $q^{1/2}$, we can identify the $B_1$ term
as follows. We have, by (2.5), that
$$\Theta_2(T_{\mathbbm C}M,T_{\mathbbm C}M,\xi_{\mathbbm C})
 =\bigotimes_{n=1}^\infty
S_{q^n}(\widetilde{T_{\mathbbm{C}}M}) \otimes
\bigotimes_{m=1}^\infty \Lambda_{-q^{m-{1\over
2}}}(\widetilde{T_{\mathbbm{C}}M}-2\widetilde{\xi}_{\mathbbm{C}})
\otimes \bigotimes_{r=1}^\infty\Lambda_{q^{r-{1\over
2}}}(\widetilde{\xi}_{\mathbbm{C}})\otimes\bigotimes_{s=1}^\infty\Lambda_{q^{s}}(\widetilde{\xi}_{\mathbbm{C}})$$
$$=
\left(1-(T_{\mathbbm{C}}M-12-2\xi_{\mathbbm{C}}+4)q^{1\over2}\right)\otimes
\left(1+(\xi_{\mathbbm{C}}-2)q^{1\over2}\right)+\cdots$$
$$=
1+(-T_{\mathbbm{C}}M+3\xi_\mathbbm{C}+6)q^{1\over2}+\cdots,\eqno(2.41)$$
where the ``$\cdots$" terms are the terms involving
 $q^{j/2}$'s with $j\geq 2$.

 {}From (2.6) and (2.37)-(2.41), one finds
 $$\left\{ {\widehat{L}(TM, \nabla^{TM})\over \cosh^2({c\over 2})}\right\}^{(12)}
 =\left\{
\left[ 8\widehat{A}(TM, \nabla^{TM}){\rm ch} (
T_{\mathbbm{C}}M,\nabla^{T_{\mathbbm{C} }M})-32\widehat{A}(TM,
\nabla^{TM})\right.\right.$$
$$- 24\left.\left.\widehat{A}(TM, \nabla^{TM})\left(e^c+e^{-c}-2\right)\right]
\cosh\left({c\over 2}\right)\right\}^{(12)}. \eqno(2.42)$$

If we set $\xi={\mathbbm R}^2$ and $c=0$ in (2.42), then we get
$$\left\{ {\widehat{L}(TM, \nabla^{TM})}\right\}^{(12)}
 =\left\{
 8\widehat{A}(TM, \nabla^{TM}){\rm ch} (
T_{\mathbbm{C}}M,\nabla^{T_{\mathbbm{C} }M})-32\widehat{A}(TM,
\nabla^{TM})\right\}^{(12)} ,\eqno(2.43)$$ which is exactly the
``miraculous cancellation" formula of Alvarez-Gaum\'e and Witten
[AGW].

$$\ $$

\noindent {\bf Appendix.  A twisted cancellation formula in $8k$
dimension}
\newcommand{\C}{\mathbbm{{C}}}

$\ $

In this appendix, we present a twisted cancellation formula for
$8k$ dimensional manifolds, which can be seen as  a direct
analogue of the $8k+4$ dimensional formula (2.7). Since the
statement is parallel and the proof is almost the same, we will
only indicate the necessary modifications. In particular, we will
use the same notation as in Section 2.

 Let $M$ be an $8k$ dimensional Riemannian
manifold with Levi-Civita connection $\nabla^{TM}$. Let $V$ be a
rank $2l$ real Euclidean vector bundle over $M$ carrying with a
Euclidean connection $\nabla^V$. Let $\xi$ be a rank two real
oriented Euclidean vector bundle over $M$ carrying with a
Euclidean connection $\nabla^\xi$. Let $R^{V}=\nabla^{V,2}$ be the
curvature of $\nabla^V$ and  $c=e(\xi,\nabla^{\xi}) $ be the
Euler form associated to $(\xi,\nabla^{\xi})$.

Let $\Theta_1(T_{\mathbbm C}M,V_{\mathbbm C},\xi_{\mathbbm C})$,
$\Theta_2(T_{\mathbbm C}M,V_{\mathbbm C},\xi_{\mathbbm C})$ be
two elements defined in the same way as in (2.5) and assume they
admit Fourier expansion in the same way as in (2.6).

We can state the main result of this appendix as follows.

$ $

\noindent {\bf Theorem A.1}  {\it If the equality for the first
Pontrjagin forms $p_1(TM,\nabla^{TM})=p_1(V,\nabla^V)$ holds, then
one has the equation for $8k $-forms,
$$\left\{ {\widehat{A}(TM,\nabla^{TM}){\det}^{1/2}
\left(2\cosh\left({\sqrt{-1}\over 4\pi}R^V\right)\right) \over
\cosh^2({c\over 2})}\right\}^{(8k)} $$ $$=2^{l+2k}\sum_{r=0}^k
2^{-6r}\left\{\widehat{A}(TM,\nabla^{TM}){\rm ch}(
b_r(T_{\mathbbm C}M,V_{\mathbbm C},\xi_{\mathbbm C })
)\cosh\left({c\over 2}\right)\right\}^{(8k)}, \eqno({\rm A}.1)$$
where each $b_r(T_{\mathbbm C}M,V_{\mathbbm C},\xi_{\mathbbm C
})$, $0\leq r\leq k$, is a  canonical integral linear combination
of $B_j(T_{\mathbbm C}M,V_{\mathbbm C},\xi_{\mathbbm C })$,
$0\leq j\leq r$.}

$\ $

If $\xi={\mathbbm R}^2$ and $c=0$, (A.1) reduces to the
cancellation formula stated in [Li, p. 32].

If we take $V=TM$ and $\nabla^V=\nabla^{TM}$, we have by (A.1)
that
$$\left\{ {\widehat{L}(TM, \nabla^{TM})\over \cosh^2\left({c\over 2}\right)}\right\}^{(8k)}
=\sum_{r=0}^k2^{6k-6r}\left\{\widehat{A}(TM,\nabla^{TM}){\rm ch}(
b_r(T_{\mathbbm C}M,T_{\mathbbm C}M,\xi_{\mathbbm C }) )
\cosh\left({c\over 2}\right)\right\}^{(8k)}. \eqno({\rm A}.2)$$

Now we assume $k=1$, that is,  $\dim M=8$. In this case, by
proceeding similarly as in Section 2e), one finds,
$$\left\{ {\widehat{L}(TM, \nabla^{TM})\over \cosh^2({c\over 2})}\right\}^{(8)}=\left\{
\left[ -\widehat{A}(TM, \nabla^{TM}){\rm ch} (
T_{\mathbbm{C}}M,\nabla^{T_{\mathbbm{C} }M})+24\widehat{A}(TM,
\nabla^{TM})\right.\right.$$
$$+3\left.\left.\widehat{A}(TM, \nabla^{TM})\left(e^c+e^{-c}-2\right)\right]
\cosh\left({c\over 2}\right)\right\}^{(8)}.\eqno({\rm A}.3) $$

We leave the details of the proofs of Theorem A.1 as well as
formula (A.3) to the interested reader.

$$\ $$

\noindent {\bf \S 3.  Spin$^c$ manifolds and Rokhlin congruences
for characteristic numbers}

$\ $

In this section, we apply our twisted cancellation formula (2.7)
to Spin$^c$ manfolds to give a direct proof of the analytic
version of the Ochanine congruence [O1] stated in [LiZ, Theorem
4.2]. In fact, the result we obtain is stronger than [LiZ,
Theorem 4.2] (see (3.2) for a precise statement).

This section is organized as follows. In a), we apply Theorem 2.1
to Spin$^c$ manifolds to get a congruence formula for
characteristic numbers. In b), we recall the analytic version of
the Ochanine congruence stated in [LiZ, Theorem 4.2] and show
that it can be proved directly as a consequence of the congruence
formula stated in a).

In this section, we will use the same notations as in Section 2.

$\ $

{\bf a). A congruence formula for Spin$^c$ manifolds }

$\ $

Let $M$ be an $8k+4$ dimensional Riemannian manifold as in
Section 2. In this section, we also assume that $M$ is closed and
oriented. Moreover, we make the assumption that there is an
$8k+2$ dimensional  closed  oriented submanifold $B$ such that if
$\widetilde{c}\in H^{2}({M,{\mathbbm Z}})$ is the Poincar\'e dual
of $[B]\in H_{8k+2}({M,{\mathbbm Z}})$, then
$$\widetilde{c}\equiv w_2(TM)\ \ \ \ \ {\rm mod}\ \ 2{\mathbbm
Z},\eqno(3.1)$$ where $w_2(TM)\in H^{2}({M,{\mathbbm Z}_2})$ is
the second Stiefel-Whitney class of $TM$.

Thus, $M$ now is a Spin$^c$ manifold. One can also show that there
exists an oriented real rank two Euclidean vector bundle $\xi$
over $M$, carrying with a Euclidean connection $\nabla^\xi$, such
that if $c=e(\xi, \nabla^\xi)$ is the Euler form associated to
$(\xi,\nabla^\xi)$, then $\widetilde{c}=[ c]$ in
$H^{2}({M,{\mathbbm Z}})$.

$\ $

\noindent {\bf Remark 3.1} If we view $\xi$ as a complex line
bundle, then $\widetilde{c}$ is the first Chern class of $\xi$.

$\ $

Now we set $V=TM$ and $\nabla^V=\nabla^{TM}$ as in Section 2b).
And for simplification we write $\Theta_2(T_{\mathbbm
C}M,\xi_{\mathbbm C})$ for $\Theta_2(T_{\mathbbm C}M,T_{\mathbbm
C}M,\xi_{\mathbbm C})$, etc. Then $\Theta_2(T_{\mathbbm
C}M,{\mathbbm C}^2)$ is exactly  the (complexification of)
$\Theta_2(TM)$ in [LiZ].

Let $B\cdot B$ be the self-intersection of $B$ in $M$. It can be
thought of as an $8k$ dimensional closed oriented manifold.

We can now state the main result of this section as follows.

$\ $

\noindent {\bf Theorem 3.2} {\it The following congruence formula
holds, $${{\rm Sign}(M)-{\rm Sign}(B\cdot B)\over 8}\equiv \int_M
\widehat{A}(TM,\nabla^{TM}){\rm ch}\left( b_k(T_{\mathbbm
C}M,\xi_{\mathbbm C }) \right) \cosh\left({c\over 2}\right)\ \ \
\ {\rm mod}\ \ 64{\mathbbm Z} ,\eqno(3.2)$$ where $b_k(T_{\mathbbm
C}M,\xi_{\mathbbm C })$ is the same term appearing in the right
hand side of $(2.40)$, and can be canonically expressed through
an integral linear combination of $B_j(T_{\mathbbm
C}M,\xi_{\mathbbm C})$, $0\leq j\leq k$.}

$\ $

 {\it Proof.} Since $[c]\in H^{2}({M,{\mathbbm Z}})$ is the Poincar\'e dual of $[B]\in
H_{8k+2}({M,{\mathbbm Z}})$, a direct computation shows that
([O1])
$$\int_M{\widehat{L}(TM, \nabla^{TM})\over \cosh^2({c\over 2})}={\rm Sign}(M)-{\rm Sign}(B\cdot
B).\eqno(3.3)$$

Formula (3.2) follows directly from (2.40), (3.3) and the
integrality result of Atiyah and Hirzebruch [AtH] stating that in
the current Spin$^c$ situation, one has
$$\int_M
\widehat{A}(TM,\nabla^{TM}){\rm ch}\left( b_j(T_{\mathbbm
C}M,\xi_{\mathbbm C }) \right) \cosh\left({c\over 2}\right)\in
{\mathbbm Z}\eqno(3.4)$$ for any $0\leq j\leq k$.  $\  $ Q.E.D.

$\ $

{\bf b). A proof of the Ochanine congruence formula}

$\ $

We need only to prove the analytic version of the Ochanine
congruence [LiZ], an equivalent version of which is stated as
[LiZ, Theorem 4.2]. So we first recall the statement of [LiZ,
Theorem 4.2] in our notation.

$\ $

\noindent {\bf Theorem 3.3} (Liu and Zhang [LiZ, Theorem 4.2])
{\it The following congruence  formula holds, $${{\rm
Sign}(M)-{\rm Sign}(B\cdot B)\over 8}\equiv \int_M
\widehat{A}(TM,\nabla^{TM}){\rm ch}\left( b_k(T_{\mathbbm
C}M+{\mathbbm C}^2-\xi_{\mathbbm C},{\mathbbm C}^2) \right)
\cosh\left({c\over 2}\right)\  {\rm mod}\ 2{\mathbbm Z}
.\eqno(3.5)$$}

{\it Proof.} By (3.2), one need only to prove that
$$\int_M
\widehat{A}(TM,\nabla^{TM}){\rm ch}\left( b_k(T_{\mathbbm
C}M,\xi_{\mathbbm C }) \right) \cosh\left({c\over 2}\right)$$
$$-\int_M \widehat{A}(TM,\nabla^{TM}){\rm ch}\left(
b_k(T_{\mathbbm C}M+{\mathbbm C}^2-\xi_{\mathbbm C},{\mathbbm
C}^2) \right) \cosh\left({c\over 2}\right)\in 2{\mathbbm
Z}.\eqno(3.6)$$

To prove (3.6), we first compare $\Theta_2(T_{\mathbbm
C}M,\xi_{\mathbbm C })$ and $\Theta_2(T_{\mathbbm C}M+{\mathbbm
C}^2-\xi_{\mathbbm C},{\mathbbm C}^2)$.

By (2.3) and (2.5), one deduces that
$$\Theta_2(T_{\mathbbm C}M+{\mathbbm
C}^2-\xi_{\mathbbm C},{\mathbbm C}^2)=\bigotimes_{n=1}^\infty
S_{q^n}(\widetilde{T_{\mathbbm{C}}M}-\widetilde{\xi}_{\mathbbm{C}})
\otimes \bigotimes_{m=1}^\infty \Lambda_{-q^{m-{1\over
2}}}(\widetilde{T_{\mathbbm{C}}M}-\widetilde{\xi}_{\mathbbm{C}})$$
$$=\Theta_2(T_{\mathbbm C}M,{\mathbbm C}^2)\otimes
{\bigotimes_{n=1}^{\infty}\Lambda_{-q^n}(\widetilde{\xi}_{\mathbbm
C })\over \bigotimes_{m=1}^{\infty}\Lambda_{-q^{m-{1\over
2}}}(\widetilde{\xi}_{\mathbbm C })},\eqno(3.7)$$
$$\Theta_2(T_{\mathbbm C}M,\xi_{\mathbbm C })=\Theta_2(T_{\mathbbm C}M,{\mathbbm C}^2)\otimes
{\bigotimes_{r=1}^{\infty}\Lambda_{q^{r-{1\over
2}}}(\widetilde{\xi}_{\mathbbm C })\otimes
\bigotimes_{s=1}^{\infty}\Lambda_{q^{s}}(\widetilde{\xi}_{\mathbbm
C })\over \left( \bigotimes_{m=1}^{\infty}\Lambda_{-q^{m-{1\over
2}}}(\widetilde{\xi}_{\mathbbm C }) \right) ^2}.\eqno(3.8)$$

{}From (3.7) and (3.8), one finds,
$$\Theta_2(T_{\mathbbm C}M,\xi_{\mathbbm C })=\Theta_2(T_{\mathbbm C}M+{\mathbbm
C}^2-\xi_{\mathbbm C},{\mathbbm
C}^2)\otimes{\bigotimes_{r=1}^{\infty}\Lambda_{q^{r-{1\over
2}}}(\widetilde{\xi}_{\mathbbm C })\otimes
\bigotimes_{s=1}^{\infty}\Lambda_{q^{s}}(\widetilde{\xi}_{\mathbbm
C })\over \bigotimes_{m=1}^{\infty}\Lambda_{-q^{m-{1\over
2}}}(\widetilde{\xi}_{\mathbbm C
})\otimes\bigotimes_{n=1}^{\infty}\Lambda_{-q^n}(\widetilde{\xi}_{\mathbbm
C })}.\eqno(3.9)$$

{} From (2.3) and the fact that ${\rm rk}(\xi)=2$,  one verifies
directly that  for any integer $r\geq 1$,
$$\Lambda_{q^r}(\widetilde{\xi}_{\mathbbm C}) \equiv \Lambda_{-q^r}(\widetilde{\xi}_{\mathbbm C})
\ \ \ \ {\rm mod}\ \ 2q^r\widetilde{\xi}_{\mathbbm C}{\mathbbm
Z}[[q^{r}]],$$
$$\Lambda_{q^{r-{1\over 2}}}(\widetilde{\xi}_{\mathbbm C}) \equiv
\Lambda_{-q^{r-{1\over 2}}}(\widetilde{\xi}_{\mathbbm C}) \ \ \ \
{\rm mod}\ \ 2q^{r-{1\over 2}}\widetilde{\xi}_{\mathbbm
C}{\mathbbm Z}[[q^{r-{1\over 2}}]].\eqno(3.10)$$

Let ${\mathbbm Z}[T_{\mathbbm C}M,{\xi}_{\mathbbm C}]$ denote the
ring with integral coefficients generated by the exterior  as
well as symmetric powers of $T_{\mathbbm C}M$ and ${\xi}_{\mathbbm
C}$.

{}From (3.9) and (3.10), one obtains,
$$\Theta_2(T_{\mathbbm C}M,\xi_{\mathbbm C })\equiv \Theta_2(T_{\mathbbm C}M+{\mathbbm
C}^2-\xi_{\mathbbm C},{\mathbbm C}^2)\ \ \ \ {\rm mod}\ \ \
2q^{1\over 2}\widetilde{\xi}_{\mathbbm C}{\mathbbm Z}[T_{\mathbbm
C}M,{\xi}_{\mathbbm C}][[q^{1\over 2}]].\eqno(3.11)$$

Thus, if one expands $\Theta_2(T_{\mathbbm C}M+{\mathbbm
C}^2-\xi_{\mathbbm C},{\mathbbm C}^2)$ as
$$\Theta_2(T_{\mathbbm C}M+{\mathbbm
C}^2-\xi_{\mathbbm C},{\mathbbm C}^2)=B_0(T_{\mathbbm
C}M+{\mathbbm C}^2-\xi_{\mathbbm C},{\mathbbm
C}^2)+B_1(T_{\mathbbm C}M+{\mathbbm C}^2-\xi_{\mathbbm
C},{\mathbbm C}^2)q^{1\over 2}+\cdots,\eqno(3.12)$$ one gets that
for any $j\geq 1$,
$$B_j(T_{\mathbbm C}M,\xi_{\mathbbm C})\equiv B_j(T_{\mathbbm C}M+{\mathbbm
C}^2-\xi_{\mathbbm C},{\mathbbm C}^2)\ \ \ \ {\rm mod}\ \ \
2\widetilde{\xi}_{\mathbbm C}{\mathbbm Z}[T_{\mathbbm
C}M,{\xi}_{\mathbbm C}].\eqno(3.13)$$

By (3.13) and by proceeding the induction arguments as in the
proof of Theorem 2.1 (Compare with [L]), one then sees easily
that for any integer $r$ such that $0\leq r\leq k$, one has
$$b_r(T_{\mathbbm C}M,\xi_{\mathbbm C})=b_r(T_{\mathbbm C}M+{\mathbbm
C}^2-\xi_{\mathbbm C},{\mathbbm C}^2)+2\widetilde{\xi}_{\mathbbm
C}C_r,\eqno(3.14)$$ for some $C_r\in{\mathbbm Z}[T_{\mathbbm
C}M,{\xi}_{\mathbbm C}]$.

Formula (3.6) then follows from (3.2), (3.14) and the
Atiyah-Hirzebruch integrality [AtH] stating that
$\int_M\widehat{A}(TM,\nabla^{TM}){\rm
ch}(\widetilde{\xi}_{\mathbbm C}C_r)\cosh({c\over 2})\in
{\mathbbm Z}$ {for any}  $0\leq r\leq k$.

The proof of Theorem 3.3 is completed. $\ $ Q.E.D.

$\ $

\noindent {\bf Remark 3.4} The factor $\widetilde{\xi}_{\mathbbm
C}$ appearing in (3.10) and (3.14) will play an important role
 in the discussion of the Finashin congruence formula in the next section.

$$\ $$

\noindent {\bf \S 4. $\eta$ invariants and a proof of the Finashin
congruence formula}

$\ $

 In this section, we combine the results in Sections 2 and 3 with the
results in [LiZ] and [Z2] to give a direct analytic proof of the
Finashin congruence formula [F].

This section is organized as follows. In a), we recall the
original statement of the Finashin congruence [F] as well as an
equivalent analytic version given in [LiZ, Theorem 4.1]. In b),
we recall a Rokhlin type congruence formula from [Z2, Theorem
3.2] in the form which is useful for the current situation. In
c), we establish a cancellation formula for certain
characteristic numbers on $8k+2$ dimensional manifolds. This
cancellation formula will be used in d), where we complete the
proof of the  analytic version of the Finashin congruence
recalled in a).

$\ $

{\bf a). The Finashin congruence and an analytic version of it}

$\ $

Let $M$ be an $8k+4$ dimensional closed smooth oriented manifold.
Let $B$ be an $8k+2$ dimensional closed smooth submanifold in $M$
such that $[B]\in H_{8k+2}(M,{\mathbbm Z}_2)$ is Poincar\'e dual
to $w_2(TM)\in H^2(M,{\mathbbm Z}_2)$. In this section, we make
the assumption that $B$ is nonorientable, as the case where $B$ is
orientable has been discussed in Section 3.

Under the above assumptions, $M\setminus B$ is oriented and spin.
We fix a spin structure on $M\setminus B$. Then it induces
canonically a pin$^-$ structure on $B$ (cf. [KT, Lemma 6.2]).

Let $o(TB)$ be the orientation bundle of $TB$. Let $L$ be the
rank two real vector bundle over $B$ defined by $L=o(TB)+{\mathbbm
R }$ (here ``$+$" stands for direct sum). Let $g^L$ be a Euclidean
metric on $L$ and denote by $L_1=\{ l\in L\,: \| l\|\leq 1 \}$
the associated unit disc bundle. Then $\partial L_1$ carries a
canonically induced spin structure from the pin$^-$ structure on
$B$ (cf. [F] and [KT]).

As an $8k+3$ dimensional oriented spin manifold, $-\partial L_1$
bounds an $8k+4$ dimensional oriented spin manifold $Z$. Following
[F], one defines
$$\Phi(B)\equiv {{\rm Sign} (Z)\over 8}\ \ \ \ {\rm mod}\ \ 2{\mathbbm
Z}.\eqno(4.1)$$ It is clear from the Ochanine divisibility (cf.
[O]) that $\Phi(B)$ is well-defined. In [F], Finashin shows that
it is a pin$^-$ cobordism invariant of $B$.

Let $B\cdot B$ be the self-intersection of $B$ in $M$. Then
$B\cdot B$ can be thought of as an $8k$ dimensional closed
oriented manifold (cf. [F], see also Section 4d)).

We can now state the original Finashin congruence as follows.

$\ $

\noindent {\bf Theorem 4.1} (Finashin [F]) {\it The following
congruence formula holds, $$ {{\rm Sign }(M)-{\rm Sign}(B\cdot
B)\over 8}\equiv \Phi(B)\ \ \ \ {\rm mod}\ \ 2{\mathbbm
Z}.\eqno(4.2)$$}

$\ $

In [LiZ], by combining the higher dimensional ``miraculous
cancellation" formula of Liu [Li] with the analytic arguments in
[Z2], Liu and Zhang give an intrinsic analytic interpretation of
the Finashin invariant $\Phi(B)$. We recall their result as
follows.

Let $g^{TB}$ be a metric on $TB$. Let $\nabla^{TB}$ be the
associated Levi-Civita connection. Let $\nabla^L$ be a Euclidean
connection on $L=o(TB)+{\mathbbm R}$.

Let $\pi:B'\rightarrow B$ be the orientable double cover of $B$.
We fix an orientation on $B'$. Then $B'$ is spin and carries an
induced spin structure from the pin$^-$ structure on $B$.

When we pull back the bundles and the associated metrics and
connections from $B$ to $B'$, we will use an extra notation ``$\,
'\, $" as indication.

Let $P:B'\rightarrow B'$ be the canonical involution on $B'$ with
respect to the double covering $\pi:B'\rightarrow B$.

Let $b_r(T_{\mathbbm C}B+o(TB)\otimes {\mathbbm C}+{\mathbbm C},
{\mathbbm C}^2)$, $0\leq r\leq k$, be the virtual Hermitian vector
bundles defined in the same way as in (2.40). They lift to
virtual vector bundles $b_r(T_{\mathbbm C}B'+{\mathbbm C}^2,
{\mathbbm C}^2)$, $0\leq r\leq k$, over $B'$, carrying with the
canonically induced $P$-invariant Hermitian connections.

Let $S(TB')$ be the bundle of spinors associated to
$(TB',g^{TB'})$. For any integer $r$ such that $0\leq r\leq k$,
let
$${D}_{B'}^{b_r(T_{\mathbbm
C}B'+{\mathbbm C}^2, {\mathbbm C}^2)}:\Gamma(S(TB')\otimes
b_r(T_{\mathbbm C}B'+{\mathbbm C}^2, {\mathbbm
C}^2))\longrightarrow \Gamma(S(TB')\otimes b_r(T_{\mathbbm
C}B'+{\mathbbm C}^2, {\mathbbm C}^2)) \eqno(4.3)$$ be the
corresponding  Dirac operator. It is a $P$-equivariant first
order elliptic differential operator and is formally
self-adjoint. As explained in [LiZ] and [Z2], ${1\over
2}(1+P){D}_{B'}^{b_r(T_{\mathbbm C}B'+{\mathbbm C}^2, {\mathbbm
C}^2)}$ determines a formally self-adjoint elliptic differential
operator $\widetilde{D}_{B}^{b_r(T_{\mathbbm C}B+o(TB)\otimes
{\mathbbm C}+{\mathbbm C}, {\mathbbm C}^2)}$ (called the twisted
Dirac operator) on $B$. When there is no confusion, we will use
the brief notation $\widetilde{D}_B^{b_r}$ to denote this twisted
Dirac operator.

Let $\overline{\eta}(\widetilde{D}_B^{b_r})$, $0\leq r\leq k$, be
the reduced $\eta$ invariant of $\widetilde{D}_B^{b_r}$ in the
sense of Atiyah, Patodi and Singer [AtPS]. By using the
Atiyah-Patodi-Singer index theorem for manifolds with boundary
[AtPS], one knows that each
$\overline{\eta}(\widetilde{D}_B^{b_r})$ is a pin$^-$ cobordism
invariant of $B$.

The following analytic interpretation of $\Phi(B)$ is proved in
[LiZ, Theorem 3.2],
$$\Phi(B)\equiv\sum_{r=0}^k
2^{6k-6r}\overline{\eta}(\widetilde{D}_B^{b_r})\ \ \ \ {\rm mod}\
\ 2{\mathbbm Z}.\eqno(4.4)$$

{}From (4.2) and (4.4), one gets the following analytic version
of the Finashin congruence formula.

$\ $

\noindent {\bf Theorem 4.2} (Liu and Zhang [LiZ, Theorem 4.1])
{\it The following congruence formula holds, $${{\rm Sign
}(M)-{\rm Sign}(B\cdot B)\over 8}\equiv \sum_{r=0}^k
2^{6k-6r}\overline{\eta}(\widetilde{D}_B^{b_r})\ \ \ \ {\rm mod}\
\ 2{\mathbbm Z}.\eqno(4.5)$$}

$\ $

In the rest of this section, we will give a direct proof of (4.5)
by combining the twisted cancellation formula (2.7) with the
analytic Rokhlin congruence formula for $KO$ characteristic
numbers proved in [Z2], which we recall in the next subsection.

$\ $

\noindent {\bf b). Rokhlin congruences for $KO$ characteristic
numbers associated to $b_r$, $0\leq r\leq k$}

$\ $

We make the same assumptions and use the same notations as before.

Let $N\rightarrow B$ be the normal bundle to $B$ in $M$. Let
$g^{TM}$ be a metric on $TM$, with the associated Levi-Civita
connection denoted by $\nabla^{TM}$. Without loss of generality
we assume that the following orthogonal splitting on $B$ holds,
$$TM|_B=TB\oplus N,\ \ \ \ \ \ \ g^{TM}=g^{TB}\oplus
g^N,\eqno(4.6)$$ where $g^N$ is the induced metric on $N$. Let
$\nabla^N$ be the Euclidean connection on $N$ induced from
$\nabla^{TM}|_B$.

Since $M$ is oriented, one has the equality for the first
Stiefel-Whitney classes,
$$w_1(TB)=w_1(N).\eqno(4.7)$$
In particular, if $o(N)$ is the orientation bundle of $N$, then
$o(N)=o(TB)$.

Let $N'=\pi^*N$ be the pull back of $N\rightarrow B$ to $B'$.
Then $N'$ is a rank two real orientable vector bundle over $B'$,
carrying with a pull back Euclidean structure, as well as a pull
back Euclidean connection $\nabla^{N'}$. Since we have fixed an
orientation on $TB'$ in the previous subsection, we see that $N'$
carries an induced orientation from the pull back of the splitting
(4.6) and from the orientation on $\pi^*(TM|_B)$.

Let $e'=e(N',\nabla^{N'})\in \Omega^2(B')$ be the Euler form
associated to $(N',\nabla^{N'})$. Then ${1\over 2}(1+P)e'$
determines an element $e\in \Omega^2(B)\otimes o(TB)$, which is
the Euler form associated to $(N,\nabla^N)$.

The following Rokhlin type congruence formula for $KO$
characteristic numbers associated to $b_r(T_{\mathbbm
C}M,{\mathbbm C}^2)$'s is a direct consequence of the general
congruence formula proved in [Z2, Section 3].

$\ $

\noindent {\bf Theorem 4.3} (Zhang [Z2, Theorem 3.2]) {\it For any
integer $r$ such that $0\leq r\leq k$, the following congruence
formula holds,
$$\int_M\widehat{A}(TM,\nabla^{TM}){\rm ch}(b_r(T_{\mathbbm
C}M,{\mathbbm C}^2))\equiv
\overline{\eta}(\widetilde{D}_B^{b_r})$$
$$+\int_B\widehat{A}(TB,\nabla^{TB}){ {\rm ch}(b_r(T_{\mathbbm C}B+N_{\mathbbm C},{\mathbbm C}^2) )-
\cosh({e\over 2}){\rm ch}(b_r(T_{\mathbbm C}B+o(TB)\otimes
{\mathbbm C}+{\mathbbm C},{\mathbbm C}^2) )\over 2\sinh({e\over
2})}\ \ \ {\rm mod}\ \ 2{\mathbbm Z}.\eqno(4.8)$$}

$\ $

\noindent {\bf Remark 4.4} The more precise expression of the
integration in the right hand side of (4.8) is
$${1\over 2}\int_{B'}\widehat{A}(TB',\nabla^{TB'}){ {\rm ch}(b_r(T_{\mathbbm C}B'+N_{\mathbbm C}',{\mathbbm C}^2) )-
\cosh({e'\over 2}){\rm ch}(b_r(T_{\mathbbm C}B' +{\mathbbm
C}^2,{\mathbbm C}^2) )\over 2\sinh({e'\over 2})}.\eqno(4.9)$$

$\ $

\noindent {\bf Remark 4.5} The proof of [Z2, Theorem 3.2] in [Z2]
uses in an essential way the techniques developed by
Bismut-Cheeger [BC] and Dai [D] on the computation of the
adiabatic limits of $\eta$ invariants of Dirac operators.

$\ $

{\bf c). A cancellation formula for characteristic numbers on
$8k+2$ dimensional manifolds}

$\ $

We assume for a moment that $B$ is oriented so that we are in the
situation discussed in Section 3.

We first set $\xi={\mathbbm R}^2$ and $c=0$ in (2.40) to get
$$\left\{ {\widehat{L}(TM, \nabla^{TM})}\right\}^{(8k+4)}
=8\sum_{r=0}^k2^{6k-6r}\left\{ \widehat{A}(TM,\nabla^{TM}){\rm
ch}( b_r(T_{\mathbbm C}M,{\mathbbm C }^2) )
\right\}^{(8k+4)},\eqno(4.10)$$ which was first proved in [Li].

{}From (2.40) and (4.10), one gets
$${1\over 8}\int_M\left( 1-{1\over\cosh^2\left({c\over
2}\right)}\right)\widehat{L}(TM,\nabla^{TM})$$
$$=\sum_{r=0}^k 2^{6k-6r}\int_M\widehat{A}(TM,\nabla^{TM})\left( {\rm
ch}( b_r(T_{\mathbbm C}M,{\mathbbm C }^2)) -\cosh\left({c\over
2}\right){\rm ch}\left( b_r(T_{\mathbbm C}M,\xi_{\mathbbm C })
\right)\right).\eqno(4.11)$$

Now since $[c]$ is Poincar\'e dual to $[B]\in
H_{8k+2}(M,{\mathbbm Z})$, one sees that if $i:B\hookrightarrow
M$ denotes the canonical embedding, then $i^*\xi=N$ and
$i^*[c]=[e]$, the Euler class of $N$.

{}From (2.1), [Hi, (9.3)] and the Chern-Weil theorem (cf. [Z3,
Chap. 1]), one deduces that
$${1\over 8}\int_M\left( 1-{1\over\cosh^2\left({c\over
2}\right)}\right)\widehat{L}(TM,\nabla^{TM})={1\over
8}\int_B\widehat{L}(TB,\nabla^{TB}){e\over \tanh\left(e\over 2
\right)}{\sinh^2\left({e\over 2}\right)\over \cosh^2\left({e\over
2}\right)}{1\over e}$$
$$={1\over
8}\int_B\widehat{L}(TB,\nabla^{TB}){\sinh\left({e\over
2}\right)\over \cosh\left({e\over 2}\right)},\eqno(4.12)$$ and
that for any integer $r$ such that $0\leq r\leq k$,
$$\int_M\widehat{A}(TM,\nabla^{TM})\left( {\rm
ch}( b_r(T_{\mathbbm C}M,{\mathbbm C }^2)) -\cosh\left({c\over
2}\right){\rm ch}( b_r(T_{\mathbbm C}M,\xi_{\mathbbm C })
)\right)$$
$$=\int_B\widehat{A}(TB,\nabla^{TB}){{\rm
ch}\left( b_r(T_{\mathbbm C}B+N_{\mathbbm C},{\mathbbm C
}^2)\right) -\cosh\left({e\over 2}\right){\rm ch}\left(
b_r(T_{\mathbbm C}B+N_{\mathbbm C},N_{\mathbbm C })\right)\over
2\sinh\left({e\over 2}\right) }.\eqno(4.13)$$

{}From (4.11)-(4.13), one gets,
$${1\over
8}\int_B\widehat{L}(TB,\nabla^{TB}){\sinh\left({e\over
2}\right)\over \cosh\left({e\over 2}\right)}$$ $$=\sum_{r=0}^k
2^{6k-6r}\int_B\widehat{A}(TB,\nabla^{TB}){{\rm ch}\left(
b_r(T_{\mathbbm C}B+N_{\mathbbm C},{\mathbbm C }^2)\right)
-\cosh\left({e\over 2}\right){\rm ch}\left( b_r(T_{\mathbbm
C}B+N_{\mathbbm C},N_{\mathbbm C })\right)\over
2\sinh\left({e\over 2}\right) }.\eqno(4.14)$$

Now we claim that (4.14) holds for any $8k+2$ dimensional closed
oriented Riemannian manifold $B$ and a rank two real oriented
Euclidean vector bundle $N$ over $B$.

In fact, for such a pair $(B,N)$, one can always take the unit
disc bundle $N_1$ of $N$. Then $\partial N_1$ is an $8k+3$
dimensional oriented spin manifold and $-\partial N_1$ bounds an
$8k+4$ dimensional oriented spin manifold $Z$. One can take
$M=N_1\cup_{\partial N_1} Z$ to get (4.14).

In the next subsection, we will apply (4.14) to the pair
$(B',N')$ discussed in Section 4b).

$\ $

\noindent {\bf Remark 4.6} Formula (4.14) actually holds on the
level of differential forms, and one can prove this directly by
using the modular invariance method, without passing to the
cobordism argument. We leave this to the interested reader.

$\ $

{\bf d). A proof of Theorem 4.2}

$\ $

We now come back to the situation of the subsections  4a) and 4b).

Recall that $B\cdot B$ denotes the self-intersection of $B$ in
$M$. It can be constructed as follows: take a transversal section
$X$ of $N$, then
$$B\cdot B=\{ b\in B\,:\,X(b)=0 \}.\eqno(4.15)$$
Let  $X'=\pi^*X$ be the pull back of $X$ over $B'$. Then $X'$ is
a transversal section of $N'$ and
$$ B'\cdot B'=\{ b'\in B'\,:\,X'(b')=0 \}\eqno(4.16)$$ is a double
cover of $B\cdot B$. Let $N'_{B'\cdot B'}$ be the normal bundle to
$B'\cdot B'$ in $B'$, then $N'_{B'\cdot B'}=N'|_{B'\cdot B'}$.
Thus, $B'\cdot B'$ carries a canonically induced orientation from
those of $TB'|_{B'\cdot B'}$ and $N'|_{B'\cdot B'}$. Moreover,
the canonical involution $P$ preserves the orientation on
$B'\cdot B'$ and thus induces an orientation on $B\cdot B$.

By a simple application of the Hirzebruch Signature theorem (cf.
[Hi]) to $B\cdot B$ and $B'\cdot B'$, one gets
$${\rm Sign}(B'\cdot B')=2\,{\rm Sign} (B\cdot B).\eqno(4.17)$$

Let $e(N'_{B'\cdot B'})$ denote the Euler class of $N'_{B'\cdot
B'}$.

{}From (2.1), [Hi, (9.3)], the Chern-Weil theorem (cf. [Z3, Chap.
1]) and the Hirzebruch Signature theorem (cf. [Hi]), one deduces
that
$$\int_{B'}\widehat{L}(TB',\nabla^{TB'}){\sinh\left({e'\over
2}\right)\over \cosh\left({e'\over
2}\right)}=\left\langle{\widehat{L}(T(B'\cdot B') )\over
e(N'_{B'\cdot B'})} {e(N'_{B'\cdot B'})\over \tanh\left(
{e(N'_{B'\cdot B'})\over 2}\right)}{\sinh\left({e(N'_{B'\cdot
B'})\over 2}\right)\over \cosh\left({e(N'_{B'\cdot B'})\over
2}\right)},[B'\cdot B']\right\rangle$$
$$=\left\langle\widehat{L}(T(B'\cdot B') ),[B'\cdot B']\right\rangle
={\rm Sign}(B'\cdot B').\eqno(4.18)$$

{}From (4.8)-(4.10), (4.17), (4.18) and (4.14) (when applied to
the pair $(B',N')$), one deduces that when mod 2${\mathbbm Z}$,
one has
$${{\rm Sign}(M)-{\rm Sign}(B\cdot B)\over 8}-\sum_{r=0}^k
2^{6k-6r}\overline{\eta}(\widetilde{D}_B^{b_r})=-{1\over
16}\int_{B'}\widehat{L}(TB',\nabla^{TB'}){\sinh\left({e'\over
2}\right)\over \cosh\left({e'\over 2}\right)}$$
$$+{1\over 2}\sum_{r=0}^k
2^{6k-6r}\int_{B'}\widehat{A}(TB',\nabla^{TB'}){{\rm ch}\left(
b_r(T_{\mathbbm C}B'+N'_{\mathbbm C},{\mathbbm C }^2)\right)
-\cosh\left({e'\over 2}\right){\rm ch}\left( b_r(T_{\mathbbm
C}B'+{\mathbbm C}^2,{\mathbbm C }^2)\right)\over
2\sinh\left({e'\over 2}\right) }$$
$$={1\over 2}\sum_{r=0}^k
2^{6k-6r}\int_{B'}\widehat{A}(TB',\nabla^{TB'}){{\rm ch}\left(
b_r(T_{\mathbbm C}B'+N'_{\mathbbm C},{\mathbbm C }^2)\right)
-\cosh\left({e'\over 2}\right){\rm ch}\left( b_r(T_{\mathbbm
C}B'+{\mathbbm C}^2,{\mathbbm C }^2)\right)\over
2\sinh\left({e'\over 2}\right) }$$ $$-{1\over 2}\sum_{r=0}^k
2^{6k-6r}\int_{B'}\widehat{A}(TB',\nabla^{TB'}){{\rm ch}\left(
b_r(T_{\mathbbm C}B'+N'_{\mathbbm C},{\mathbbm C }^2)\right)
-\cosh\left({e'\over 2}\right){\rm ch}\left( b_r(T_{\mathbbm
C}B'+N'_{\mathbbm C},N'_{\mathbbm C })\right)\over
2\sinh\left({e'\over 2}\right) }$$
$$={1\over 2}\int_{B'}\widehat{A}(TB',\nabla^{TB'})
{\cosh\left({e'\over 2}\right)\over 2\sinh\left({e'\over
2}\right)} \left({\rm ch}(b_r(T_{\mathbbm C}B'+N'_{\mathbbm C
},N'_{\mathbbm  C}))- {\rm ch}(b_r(T_{\mathbbm C}B'+{\mathbbm C
}^2,{\mathbbm  C}^2))\right).\eqno(4.19)$$ {}From (4.19), one sees
that in order to prove Theorem 4.2, one need only to prove the
following result.

$\ $

\noindent {\bf Lemma 4.7} {\it For any integer $r$ such that
$0\leq r\leq k$, one has $${1\over
2}\int_{B'}\widehat{A}(TB',\nabla^{TB'}) {\cosh\left({e'\over
2}\right)\over 2\sinh\left({e'\over 2}\right)} \left({\rm
ch}(b_r(T_{\mathbbm C}B'+N'_{\mathbbm C },N'_{\mathbbm  C}))-
{\rm ch}(b_r(T_{\mathbbm C}B'+{\mathbbm C }^2,{\mathbbm
C}^2))\right)\in 2{\mathbbm Z}.\eqno(4.20)$$}

{\it Proof.} Let ${\mathbbm Z}[T_{\mathbbm C}B',N'_{\mathbbm C}]$
be the ring with integral coefficients generated by the exterior
and symmetric powers of $T_{\mathbbm C}B'$ and $N'_{\mathbbm C}$.

By an obvious analogue of (3.14), one has that for any integer
$r$ such that $0\leq r\leq k$,
$$b_r(T_{\mathbbm C}B'+N'_{\mathbbm C
},N'_{\mathbbm  C})- b_r(T_{\mathbbm C}B'+{\mathbbm C
}^2,{\mathbbm  C}^2)=2\widetilde{N}'_{\mathbbm
C}T'_r,\eqno(4.21)$$ for some $T'_r\in {\mathbbm Z}[T_{\mathbbm
C}B',N'_{\mathbbm C}]$.

{}From (4.21), one deduces that
$${\cosh\left({e'\over 2}\right)\over 2\sinh\left({e'\over
2}\right)} \left({\rm ch}(b_r(T_{\mathbbm C}B'+N'_{\mathbbm C
},N'_{\mathbbm  C}))- {\rm ch}(b_r(T_{\mathbbm C}B'+{\mathbbm C
}^2,{\mathbbm  C}^2))\right)$$
$$={\cosh\left({e'\over 2}\right)\over \sinh\left({e'\over
2}\right)}\left(\exp(e')+\exp(-e')-2\right){\rm ch}(T'_r)
={\cosh\left({e'\over 2}\right)\over \sinh\left({e'\over
2}\right)}4\sinh^2\left({e'\over 2}\right){\rm ch}(T'_r)$$ $$
=2\sinh(e'){\rm ch}(T'_r).\eqno(4.22)$$

{}From (4.22), [Hi, (9.3)] and the Chern-Weil theorem (cf. [Z3,
Chap. 1]), one deduces that
$${1\over 2}\int_{B'}\widehat{A}(TB',\nabla^{TB'})
{\cosh\left({e'\over 2}\right)\over 2\sinh\left({e'\over
2}\right)} \left({\rm ch}(b_r(T_{\mathbbm C}B'+N'_{\mathbbm C
},N'_{\mathbbm  C}))- {\rm ch}(b_r(T_{\mathbbm C}B'+{\mathbbm C
}^2,{\mathbbm  C}^2))\right)$$
$$=\int_{B'}\widehat{A}(TB',\nabla^{TB'}){\rm ch}(T'_r)\sinh(e')$$
$$=\left\langle \widehat{A}(T(B'\cdot B')){\rm ch}(i^*_{B'\cdot B'}T'_r)
{\sinh (e(N'_{B'\cdot B'}))\over e(N'_{B'\cdot
B'})}{e(N'_{B'\cdot B'})\over 2\sinh\left({e(N'_{B'\cdot
B'})\over 2}\right)},[B'\cdot B']\right\rangle$$
$$=\left\langle \widehat{A}(T(B'\cdot B')){\rm ch}(i^*_{B'\cdot
B'}T'_r)\cosh\left({e(N'_{B'\cdot B'})\over 2}\right),[B'\cdot
B']\right\rangle,\eqno(4.23)$$ where $i^*_{B'\cdot B'}:B'\cdot
B'\hookrightarrow B'$ denotes the canonical embedding.

It is clear that $i^*_{B'\cdot B'}T'_r$ is invariant under the
canonical involution $P$ and thus induces a virtual complex
vector bundle over $B\cdot B$. We denote it by $T_r(B\cdot B)$.

Let $N_{B\cdot B}$ be the normal bundle to $B\cdot B$ in $B$.
Then $e(N'_{B'\cdot B'})$ is the pull back of the Euler class
$e(N_{B\cdot B})$ of $N_{B\cdot B}$ through the covering map
$B'\cdot B'\rightarrow B\cdot B$. Moreover, it is clear that the
total Pontrjagin class of the oriented real vector bundle
$N_{B\cdot B}\oplus o(N_{B\cdot B})$ is given by
$$p(N_{B\cdot B}\oplus o(N_{B\cdot B}))=1+\left(e(N_{B\cdot
B})\right)^2.\eqno(4.24)$$

With these notations, one gets
$$\left\langle \widehat{A}(T(B'\cdot B')){\rm ch}(i^*_{B'\cdot
B'}T'_r)\cosh\left({e(N'_{B'\cdot B'})\over 2}\right),[B'\cdot
B']\right\rangle$$ $$=2\left\langle \widehat{A}(T(B\cdot B)){\rm
ch}(T_r(B\cdot B))\cosh\left({e(N_{B\cdot B})\over
2}\right),[B\cdot B]\right\rangle.\eqno(4.25)$$

Now as $B$ is pin$^-$, one has the following equality for the
Stiefel-Whitney classes (cf. [KT] and [F]),
$$w_2(TB)+(w_1(TB))^2=0.\eqno(4.26)$$
By (4.7) and (4.26), one gets
$$w_2(TB)+(w_1(N))^2=0.\eqno(4.27)$$
Pulling back (4.27) to $B\cdot B$, one gets
$$w_2(T(B\cdot B))+w_2(N_{B\cdot B})+(w_1(N_{B\cdot
B}))^2=0.\eqno(4.28)$$ On the other hand, one has
$$w_2(N_{B\cdot B}\oplus o(N_{B\cdot B}))=w_2(N_{B\cdot B})+(w_1(N_{B\cdot
B}))^2.\eqno(4.29)$$

{}From (4.28) and (4.29), one finds
$$w_2(T(B\cdot B))=w_2(N_{B\cdot B}\oplus o(N_{B\cdot
B})).\eqno(4.30)$$ From (4.24), (4.30) and the obvious fact that
$T_r(B\cdot B)$ is indeed the complexification of some (virtual)
real vector bundle over $B\cdot B$, one  then applies a result of
Mayer [M, Satz 3.2(vi)] to conclude that
$$2\left\langle \widehat{A}(T(B\cdot B)){\rm
ch}(T_r(B\cdot B))\cosh\left({e(N_{B\cdot B})\over
2}\right),[B\cdot B]\right\rangle\in 2{\mathbbm Z}.\eqno(4.31)$$

Formula (4.20) then follows from (4.23), (4.25) and (4.31).

The proof of Lemma 4.7 is completed. \ \ \  Q.E.D.

The proof of Theorem 4.2 is thus also completed.\ \ \ Q.E.D.

$\ $

\noindent {\bf Remark 4.8} It is remarkable that Mayer's result
is needed here in  the $B$ nonorientable case. The basic reason
is that $T(B\cdot B)+N_{B\cdot B}+o(N_{B\cdot B})$ is a rank
$8k+3$ real oriented  spin vector bundle. Thus the associated
spinor bundle   carries a quarternionic structure. Mayer's result
is then a direct consequence of the Atiyah-Singer index theorem
[AtS].

$$\ $$

\noindent{\bf References}

$\ $

\noindent [AGW] L. Alvarez-Gaum\'e and E. Witten, Gravitational
anomalies. {\it Nucl. Phys.} B234 (1983), 269-330.

$\ $

\noindent [At] M. F. Atiyah, {\it K-theory}. Benjamin, New York,
1967.

$\ $

\noindent [AtH] M. F. Atiyah and F. Hirzebruch, Riemann-Roch
theorems for differentiable manifolds. {\it Bull. Amer. Math.
Soc.} 65 (1959), 276-281.

$\ $

\noindent [AtPS] M. F. Atiyah, V. K. Patodi and I. M. Singer,
Spectral asymmetry and Riemannian geometry I. {\it Proc.
Cambridge Philos. Soc.} 77 (1975), 43-69.

$\ $

\noindent [AtS] M. F. Atiyah and I. M. Singer, The index of
elliptic operators III. {\it Ann. of Math.} 87 (1968), 546-604.

$\ $

\noindent [BC] J.-M. Bismut and J. Cheeger, $\eta$-invariants and
their adiabatic limits. {\it J. Amer. Math. Soc.} 2 (1989), 33-70.

$\ $

\noindent [C] K. Chandrasekharan, {\it Elliptic Functions}.
Springer-Verlag, 1985.

$\ $

\noindent [D] X. Dai, Adiabatic limits, nonmultiplicity of
signature, and Leray spectral sequence. {\it J. Amer. Math. Soc.}
4 (1991), 265-321.

$\ $

 \noindent [F] S. M. Finashin, {A Pin$^-$-cobordism
invariant and a generalization of the Rokhlin signature
congruence.} {\it Leningrad Math. J.} 2 (1991), 917-924.

$\ $



\noindent [GuM] L. Guillou and A. Marin, Une extension d'un
th\'eor\`eme de Rokhlin sur la signature. {\it C. R. Acad. Sci.
Paris, S\'erie I}. 285 (1977), 98-98.

$\ $

\noindent [HZ] F. Han and W. Zhang, Spin$^c$-manifolds and
elliptic genera.  {\it C. R. Acad. Sci. Paris, S\'erie I.}, To
appear.

$\ $

\noindent [Hi] F. Hirzebruch, {\it Topological Methods in
Algebraic Geometry}. Springer-Verlag, 1966.

$\ $

\noindent [KT] R. C. Kirby and L. R. Taylor, Pin structures on low
dimensional manifolds. in {\it Geometry of Low Dimensional
Manifolds}, Vol. 2, p. 177-242. Ed. S. K. Donaldson and C. B.
Thomas, Cambridge Univ. Press, 1990.

$\ $

\noindent [L] P. S. Landweber, Elliptic cohomology and modular
forms. in {\it Elliptic Curves and Modular Forms in Algebraic
Topology}, p. 55-68. Ed. P. S. Landweber. Lecture Notes in
Mathematics Vol. 1326, Springer-Verlag (1988).

$\ $

\noindent [Li] K. Liu, Modular invariance and characteristic
numbers. {\it Commun. Math. Phys}. 174 (1995), 29-42.

$\ $

\noindent [LiZ] K. Liu and W. Zhang, Elliptic genus and
$\eta$-invariants. {\it Inter. Math. Res. Notices} No. 8 (1994),
319-328.

$\ $

\noindent [M] K. L. Mayer, Elliptische differentialoperatoren und
ganzzahligkeitss\"atze f\"ur charakteristische zahlen. {\it
 Topology} 4 (1965), 295-313.

 $\ $

\noindent [O1] S. Ochanine, Signature modulo 16, invariants de
Kervaire g\'en\'eralis\'es et nombre caract\'eristiques dans la
$K$-th\'eorie reelle. {\it M\'emoire Soc. Math. France}, Tom. 109
(1987), 1-141.

$\ $

\noindent [O2] S. Ochanine, Elliptic genera, modular forms over
$KO_*$, and the Brown-Kervaire invariants. {\it Math. Z.} 206
(1991), 277-291.

$\ $

\noindent [R1] V. A. Rokhlin, New results in the theory of
4-dimensional manifolds. {\it Dokl. Akad. Nauk. S.S.S.R.,} 84
(1952), 221-224.

$\ $

\noindent [R2] V. A. Rokhlin, Proof of a conjecture of Gudkov.
{\it Funct. Anal. Appl.}, 6 (1972), 136-138.

$\ $



\noindent [Z1] W. Zhang, Spin$^c$-manifolds and Rokhlin
congruences.
 {\it   C. R. Acad. Sci. Paris,
S\'erie I}, 317 (1993), 689-692.

$\ $

\noindent [Z2] W. Zhang, {Circle bundles, adiabatic limits of
$\eta$ invariants and Rokhlin congruences.} {\it   Ann. Inst.
Fourier} 44 (1994), 249-270.

$\ $

\noindent [Z3] W. Zhang, {\it Lectures on Chern-Weil Theory and
Witten Deformations.} Nankai Tracks in Mathematics Vol. 4, World
Scientific, Singapore, 2001.

$$\ $$

\noindent Nankai Institute of Mathematics, Nankai University,
Tianjin 300071, PR China

$\ $

\noindent {\it E-mails:}

\noindent F. H.: hanfeiycg@yahoo.com.cn

\noindent W.Z.: weiping@nankai.edu.cn

\end{document}